\documentclass[notitlepage,11pt]{article}
\usepackage{amssymb,amsmath, cases, tikz, float, subfig}
\catcode`\@=11
\@addtoreset{equation}{section}

\catcode`\@=12
\usepackage{colortbl}%

\def\R{\mathbb{R}}

\def\f{\varphi}

\def\irn{\int_{\R^n}}

\def\f{\varphi}

\def\QED{\hfill {$\square$}\goodbreak \medskip}

\newtheorem{Theorem}{Theorem}[section]

\newtheorem{Corollary}[Theorem]{Corollary}

\newtheorem{Definition}[Theorem]{Definition}

\begin{document}

\title 
{The H\'enon-Lane-Emden system: \\
a sharp nonexistence result}

\author{Andrea Carioli\footnote{SISSA, via Bonomea, 265 -- 34136 Trieste, Italy. Email: {acarioli@sissa.it}.
{Partially supported by INDAM-GNAMPA.}}~ and
Roberta Musina\footnote{Dipartimento di Matematica e Informatica, Universit\`a di Udine,
via delle Scienze, 206 -- 33100 Udine, Italy. Email: {roberta.musina@uniud.it}. 
{Partially supported by Miur-PRIN 201274FYK7\_004.}
}}

\date{}

\maketitle

\footnotesize

\noindent

\begin{abstract}{\small \noindent
We deal with very weak positive supersolutions to the H\'enon-Lane-Emden system on neighborhoods
of the origin. In our main theorem
we prove a sharp nonexistence result.
\medskip

\medskip

}
\end{abstract}

\normalsize
\vskip1cm

\section{Introduction}
The  system of elliptic equations
\begin{equation}
\label{eq:power0}
\begin{cases}
-\Delta u= |x|^av^{p-1}\\
-\Delta v= |x|^bu^{q-1}
\end{cases}
\end{equation}
has been largely studied since Mitidieri's paper~\cite{Mit-1} appeared, in 1990.
We focus our attention on the related problem
\begin{equation*}\tag{$\mathcal P_{a,b}$}
\label{eq:HLE}
\begin{cases}
-\Delta u\ge \lambda_1|x|^av^{p-1}\\
-\Delta v\ge \lambda_2|x|^bu^{q-1}
\end{cases} 
\end{equation*}
on punctured domains $\Omega\setminus\{0\}$,
where $p,q>1$, $\Omega\subset\R^n$ is a neighborhood of the origin and $n\ge 3$. 
We are interested in nonnegative, distributional (or {\em very weak}) solutions to (\ref{eq:HLE}), accordingly with
the next definition.

\begin{Definition}
\label{D:vw}
A nontrivial and {nonnegative distributional solution} to \eqref{eq:HLE} on $\Omega\setminus\{0\}$ is a pair $u,v$ of nonnegative functions satisfying

\medskip

\centerline{$u, v\in L^1_{\rm loc}(\Omega\setminus\{0\})$, $u^{q-1}, v^{p-1}\in L^1_{\rm loc}(\Omega\setminus\{0\})$,}

\medskip
\noindent
for which there exist $\lambda_1,\lambda_2>0$ such that the inequalities in \eqref{eq:HLE} hold in the sense of distributions on 
$\Omega\setminus\{0\}$.
\end{Definition}
Problems (\ref{eq:power0}) and (\ref{eq:HLE}) change their 
nature depending on the sign of the quantity ${(p-1)(q-1)-1}$. One has to distinguish 
between the following cases:
$$
\begin{array}{lrc}
\text{(AC)}\qquad\qquad&\displaystyle{\frac{1}{p}+\frac{1}{q}}< 1&\qquad\text{[Anticoercive case]}\\
&&\\
\text{~(H)}&\displaystyle{\frac{1}{p}+\frac{1}{q}}= 1&\qquad\text{[Homogeneous case]}\\
&&\\
\text{~({C})}&\displaystyle{\frac{1}{p}+\frac{1}{q}}> 1&\qquad\text{[Coercive case]}
\end{array}
$$
In the homogeneous case (H) the parameters $\lambda_1$, $\lambda_2$
have to be regarded as (possibly nonlinear) eigenvalues and can not be {\em a priori} prescribed. If
(AC) or ({C}) applies, then one can always assume that $\lambda_1=\lambda_2=1$.

In the present paper we prove a sharp nonexistence result in the spirit
of the paper \cite{BC} by Brezis and Cabr\'e. More precisely,
for  fixed $p,q>1$ we find the region $E_{p,q}$
of parameters $a,b\in\R$, for which there exist  positive distributional solutions to  \eqref{eq:HLE}
in neighborhoods of the origin. The set $E_{p,q}$ is defined as follows.

\medskip

\noindent
$\bullet$ {\bf Anticoercive case.} If (AC) holds, then
$$E_{p,q}:=
\displaystyle{\left\{(a,b)~\big|~ a,b>-n, ~~~
\displaystyle\frac{a}{p}+\frac{b}{p'}+2> 0~,~~
\displaystyle\frac{a}{q'}+\frac{b}{q}+2> 0
\right\}}~\!.
$$

\medskip

\noindent
$\bullet$ {\bf Homogeneous case.} We put
$$E_{p,p'}:=
\displaystyle{\left\{(a,b)~\big|~ a,b>-n, ~~~
\displaystyle\frac{a}{p}+\frac{b}{p'}+2\ge 0~
\right\}}~\!.
$$

\medskip

\noindent
$\bullet$ {\bf Coercive case.} If ({C}) holds, then
$$E_{p,q}:=
\displaystyle{\left\{(a,b)~\big|~ a,b>-n, ~~~
\begin{array}{ll}
\displaystyle\frac{a+n}{p}+\frac{b+n}{q(p-1)}>n-2\\
~\\
\displaystyle\frac{a+n}{p(q-1)}+\frac{b+n}{q}>n-2
\end{array}
\right\}}~\!.
$$
We are in position to state our main result.

\begin{Theorem}
\label{T:vw}
Let $\Omega\subset \R^n$ be a  bounded domain containing the origin, and let $p,q>1$.
Then \eqref{eq:HLE} has a nontrivial and nonnegative distributional solution on $\Omega\setminus\{0\}$
if and only if $(a,b)\in {E}_{p,q}$.
\end{Theorem}

\noindent Trivially,  any distributional solution
$u\ge 0$ to 
\begin{equation}
\label{eq:BC}
-\Delta u\ge \lambda|x|^au^{p-1}
\end{equation}
gives rise to the solution $u, v=u$ to the corresponding system
\begin{equation}
\label{eq:BC_sys}
\begin{cases}
-\Delta u\ge \lambda_1|x|^av^{p-1}\\
-\Delta v\ge \lambda_2|x|^au^{p-1}~\!,
\end{cases} 
\end{equation}
for $\lambda_1=\lambda_2=\lambda$.
The converse might not be true, in general. 
Recall that by \cite[Theorem 0.1]{BC}, the inequality (\ref{eq:BC})
has no nontrivial and nonnegative distributional solutions on $\Omega\setminus\{0\}$
if $\lambda>0$, $p=3$ and $a\le -2$
(see  also \cite[Theorem~1.2]{CalM} for  $p>2$ and for more general nonlinearities). 
Thanks to Theorem \ref{T:vw}, we immediately get the following extension of \cite[Theorem 0.1]{BC}
for systems.

\begin{Corollary}
\label{C:cor}
Let $\Omega\subset \R^n$ be a  bounded domain containing the origin, and let $p>1$.
Then (\ref{eq:BC_sys}) has a nontrivial and nonnegative distributional solution on $\Omega\setminus\{0\}$
if and only if one of the following conditions is satisfied:
\begin{description}
\item$~~i)$ $p>2$ and $a>-2$
\item$~ii)$ $p=2$ and $a\ge -2$
\item$iii)$ $1<p<2$, $a>-n$ and $p<p_a:=\displaystyle \frac{2(n-1)+a}{n-2}$.
\end{description}
\end{Corollary}
Notice that $p_a$ coincides with  {\em Serrin's critical exponent} when  $a=0$.

\medskip

Theorem \ref{T:vw}, combined with the action of the  {\em Kelvin transform}
\begin{equation}
\label{eq:kelvin}
\mathcal K: L^1_{\rm loc}(\R^n\setminus\{0\})\to L^1_{\rm loc}(\R^n\setminus\{0\})~,\quad (\mathcal K w)(x)=|x|^{2-n}w\Big(\frac{x}{|x|^2}\Big)~\!,
\end{equation}
immediately leads to the following sharp Liouville-type result
(power-type solutions are computed in the appendix).

\begin{Theorem}
\label{T:B}
The system of inequalities \eqref{eq:HLE} has a positive distributional solution $u,v$ on $\R^n\setminus\{0\}$
if and only if the system (\ref{eq:power0})
has a positive power-type solution, that is, having the 
form $u(x)=c_1|x|^\alpha, v(x)=c_2|x|^\beta$.
\end{Theorem}

\medskip

Theorems \ref{T:vw} and \ref{T:B} are related to some known results. 
Serrin and Zou \cite{SZ} constructed positive radial solutions of class
$C^2(\R^n)$  under the assumptions $a=b=0$,
$n/p+n/q\le n-2$.  One can adapt the shooting method in \cite{SZ} 
to find a bounded  solution to (\ref{eq:power0})
in a ball $\Omega$ about the origin if $a,b>-2$ and $(a,b)\in E_{p,q}$.
The restriction $a,b>-2$ is  necessary to find solutions of class $C^0(\Omega)\cap C^2(\Omega\setminus\{0\})$,
see for instance \cite{BV}.  

Under the anticoercivity assumption  (AC), Bidaut-Veron and Giacomini  \cite{BVG}
investigated an equivalent Hamiltonian system to prove the existence of a radial solution
$u,v\in C^0(\R^n)\cap C^2(\R^n\setminus\{0\})$ on $\R^n$ if and only if
$a,b>-2$ and 
$$
\frac{a+n}{p}+\frac{b+n}{q}\le n-2.
$$ 
To prove the existence part in Theorem \ref{T:vw} we write
$$
E_{p,q}=E^+_{p,q}\cup E^-_{p,q}~\!,
$$
where $E^+_{p,q}$ is the set of pairs $(a,b)$ such that
\begin{equation}
\label{eq:above_line}
a,b>-n~,\qquad \frac{a+n}{p}+\frac{b+n}{q}>n-2~\!,
\end{equation}
and $E^-_{p,q}=E_{p,q}\setminus E^+_{p,q}$. If $(a,b)\in E^-_{p,q}$ then the system (\ref{eq:HLE})
admits power type solutions, see the explicit computations in the Appendix. For $(a,b)\in E^+_{p,q}$
we take a large ball $B$ about the origin containing $\Omega$ and we study the 
system
\begin{equation}
\label{eq:nuovo}
\begin{cases}
-\Delta u=\lambda_1|x|^av^{p-1}\\
-\Delta v=\lambda_2|x|^bu^{q-1} \\
u, v>0&\text{in $B$}\\
u=v=0&\text{on $\partial B$}
\end{cases} 
\end{equation}
The existence of a solution to (\ref{eq:HLE}) readily
follows from the next theorem, that might have an independent interest.  

\begin{Theorem}
\label{T:radial}
If (\ref{eq:above_line}) holds, then there exist $\lambda_1,\lambda_2>0$ such that
(\ref{eq:nuovo}) has at least a radial solution $u,v$ satisfying
\begin{equation}
\label{eq:energy}
u,v \in W^{2,1}\cap W^{1,1}_0(B)~, \qquad \int_B|x|^b|u|^q~\!dx<\infty~,\qquad \int_B|x|^a|v|^p~\!dx<\infty.
\end{equation}
\end{Theorem}

Theorem \ref{T:radial} will be proved in Section \ref{S:existence}, via variational arguments. 
A simple computation shows that $u,v$ can never be a power-type solution.
Notice that
Theorem \ref{T:radial} provides existence also for certain exponents $a, b\ngtr-2$.

\medskip

Most of the available nonexistence results 
concern the system of elliptic equations (\ref{eq:power0}) or deal with
more regular solutions. First of all, we cite the pioneering paper \cite{GS} by Gidas and Spruck,
and in particular their Theorem A.3. In \cite{BV}, Bidaut-Veron used a clever and interesting trick to prove a nonexistence result of 
classical solutions $u,v$ to (\ref{eq:power0}) on the punctured domain $\Omega\setminus\{0\}$. 
Her argument plainly covers  problem (\ref{eq:HLE}) and can be
used to prove Theorem \ref{T:B} under the additional assumption $u,v\in C^2(\R^n\setminus\{0\})$.

D'Ambrosio and Mitidieri \cite[Theorem 3.5]{D'AM} used   representation formulae
to prove the nonexistence result in Theorem \ref{T:B} for 
locally integrable distributional solutions 
on $\R^n$.
Notice however that \cite{D'AM} include
a much larger class of systems. 

The nonexistence part of Theorem \ref{T:vw} will be proved in Section \ref{S:ne}.

\bigskip

The available literature for (\ref{eq:power0}) and related problems is very extensive.
The interested reader can find exhaustive surveys in \cite{BVG, D'AM, FG}, besides remarkable results.
A number of papers (see for instance \cite{PQS, SZ2, Sou, Souto}) deal with the 
{\em H\'enon-Lane-Emden Conjecture}, that originated
from the nonexistence results in \cite{Mit-1}, \cite{Mit1}. We cite also 
\cite{BMR1, BMR2, BusMan, CDM, CFM, DPR, F, GMMY, Mit0, MS, PVVV, Ph, SZ3},
and the references therein.

\section{Proof of Theorem \ref{T:radial}: existence}
\label{S:existence}
The homogeneous case (H) is  covered by
\cite[Theorem 1.3]{CM}. 
Thus, we assume that $q\neq p'$. Since the system (\ref{eq:nuovo}) is not homogeneous,
we can fix $\lambda_1=\lambda_2=1$. To get existence, we follow the outline of the proof in \cite{CM}. 
For brevity, we will skip some details.  

Our approach is based on the formal equivalence, already noticed for instance
in \cite{ClMit}, between the system
(\ref{eq:nuovo} and 
the fourth order Navier problem
\begin{equation}
\label{eq:fourth}
\begin{cases}
-\Delta\left( |x|^{-a(p'-1)} (-\Delta u)^{p'-1} \right) =  |x|^b u^{q-1} \\
u, -\Delta u>0&\text{in $B$}\\
u=\Delta u=0&\text{on $\partial B$.}
\end{cases}
\end{equation}
We use variational methods to show that (\ref{eq:fourth}) 
admits a  radial  weak solution $u$ in a suitably defined energy space. To
conclude the proof, one only has to check
that the pair
$
u, ~ v:=|x|^{-a(p'-1)}(-\Delta u)^{p'-1}$
solves (\ref{eq:nuovo}).

The first step consists in defining
$$
W^{2,p'}_{N,\mathrm{rad}}(B;|x|^{-a(p'-1)}dx)
$$
as the completion of of the space of radial functions $u\in C^2(\overline{B})$, such that 
$$
u = 0 ~~\text{ on } \partial\Omega~,\qquad \Delta u \equiv 0~~ \text{in a neighborhood of
the origin},
$$
 with respect to the norm
$$
\|u\|=\left(\int_B|x|^{-a(p'-1)}|\Delta u|^{p'}~\!dx\right)^{\!1/p'}.
$$
We claim that the infimum
$$
m:= \inf_{\substack{u\in W^{2,p'}_{N,\text{rad}}(B;|x|^{-a(p'-1)}dx) \\ u\neq 0}}
\frac{\displaystyle\int_B|x|^{-a(p'-1)}|\Delta u|^{p'}~\!dx}
{\left(\displaystyle\int_B|x|^b|u|^q~\!dx\right)^{p'/q} }
$$
is positive and achieved by some function $ u$. 
For the sake of clarity, we distinguish the coercive case from the anticoercive one.

\medskip
\noindent
{\bf Coercive case.} If $q<p'$ we take any exponent $b_0$, such that $b_0>-n$
and
$$
(i)~~\frac{a}{p}+\frac{b_0}{p'}+2>0~,\qquad (ii)~~\frac{b_0+n}{p'}<\frac{b+n}{q}.
$$
Thanks to $(i)$, we have that  $W^{2,p'}_{N,\text{rad}}(B;|x|^{-a(p'-1)}dx)$ is compactly embedded into 
$L^{p'}(B;|x|^{b_0} dx)$ by \cite[Lemma 2.8]{CM}. On the other hand, the space $L^{p'}(B;|x|^{b_0} dx)$ is continuously embedded into
$L^{q}(B;|x|^{b} dx)$ by $(ii)$ and H\"older inequality. The claim follows by standard arguments.

\medskip
\noindent
{\bf Anticoercive case.} Fix exponents $a_0, b_0$ satisfying 
\begin{gather*}
\frac{a_0+n}{p}+\frac{b_0+n}{q}=n-2~\!,\quad {-n < a_0 \le a,\ -n < b_0 < b,}
\end{gather*}
that is possible as $(a,b)\in E_{p,q}$ and (\ref{eq:above_line}) holds.
By \cite[Theorem 4.10]{M}, we have that there exists
a constant $c>0$ such that
$$
\irn|x|^{-a_0(p'-1)}|\Delta \f|^{p'}~\!dx\ge c\left(\irn|x|^{b_0}| \f|^{q}~\!dx\right)^{p'/q}
$$
for any radially symmetric $\f\in C^\infty_c(\R^n)$. Since we are dealing with a bounded domain, it is easy to infer that
$W^{2,p'}_{N,\text{rad}}(B;|x|^{-a(p'-1)}dx)$ is compactly embedded into 
$L^q(B;|x|^b dx)$, and this proves the claim.

\bigskip

Next, let $u$ be an extremal for the infimum $m$.
Use the arguments in  \cite[Lemma 3.2]{CM} to show that 
$$
u\in W^{2,1}\cap W^{1,1}_0(B)~,\quad v:=|x|^{-a(p'-1)}|\Delta u|^{p'-2}(-\Delta u)\in W^{2,1}\cap W^{1,1}_0(B)~\!,
$$
and that, up to a Lagrange multiplier,  the pair $u,v$ is a weak solution to the system
$$
-\Delta u=|x|^a|v|^{p-2}v~,\qquad -\Delta v=|x|^b|u|^{q-2}u
$$
in the ball $B$. To check that  $u, v$ are positive on $B$
use the (standard) argument
in \cite[Lemma 3.4]{CM}. 
The proof of the existence part is complete.
\QED

\section{Proof of Theorem \ref{T:vw}: nonexistence}
\label{S:ne}
In this proof we denote by $c$ any inessential positive constant.

Let $u,v$ be a nontrivial and nonnegative distributional solution to \eqref{eq:HLE} in $\Omega\setminus\{0\}$. 
We claim that the following facts hold: 
\begin{description}
\item$~~i)$ 
$u, v\in L^1_{\rm loc}(\Omega)$, $|x|^bu^{q-1}, |x|^av^{p-1}\in L^1_{\rm loc}(\Omega)$;
\item$~ii)$ $u,v$ solve \eqref{eq:HLE} in the sense of distributions on $\Omega$;
\item$iii)$ $u,v$ are superharmonic and positive on $\Omega$;
\item$~iv)$ $a,b>-n$.
\end{description}
The first two conclusions are immediate consequences of 
\cite[Lemma 1]{BDT}. Since $u,v\in L^1_{\rm loc}(\Omega)$ solve
$-\Delta u\ge 0$, $-\Delta v\ge 0$, then $u,v$ are superharmonic by
well known and classical facts. In particular $u$ and $v$ can be assumed to be lower semicontinuous and
positive on $\Omega$, and so $iii)$ holds. Finally, since $v$ is lower semicontinuous and positive,
we can find $\delta>0$ such that 
$|x|^av^{p-1}\ge \delta |x|^a$ in a closed ball $\overline B\subset\Omega$ about the origin. Now, from 
$|x|^av^{p-1}\in L^1(B)$ we infer that the
weight $|x|^a$ is locally integrable on $B$, that is, $a>-n$. The conclusion $b>-n$
can be proved in a similar way.

Next, up to dilations we can assume that $\Omega$ contains the closure of the unit ball $B$
about the origin. Let $G^x(\cdot)$ be the Green function for $B$ and
let $h^x(\cdot)$ be its regular part, that is,
\[
G^x(y) = c_n\left[ \frac{1}{| y-x |^{n-2}} -
h^x(y)\right]~,\quad h^x(y)= \frac{|x|^{n-2}}{|y|x|^2- x|^{n-2}}~\!.
\]
We claim that
\begin{equation}\label{eq:sysIne}
\begin{gathered}
u(x) \ge \lambda_1 \int_B G^x_B(y) | y |^a v(y)^{p-1}~\!dy,\\
v(x) \ge \lambda_2 \int_B G^x_B(y) | y |^b u(y)^{q-1}~\!dy
\end{gathered}
\end{equation}
almost everywhere on $B$.
Let us prove the first of the inequalities in \eqref{eq:sysIne}, the second one being similar.
For any integer $k\ge 1$, we put
$$
f_k=\min\{\lambda_1|x|^a v^{p-1}, k\}
$$
and we introduce the unique solution $u_k$ to the problem 
$$
-\Delta u_k = f_k~, \qquad u_k  \in H^1_0(B).
$$
Green's representation formula yields 
\[
u_k(x) = \int_B G^x_B(y)f_k~\!dy,
\]
and the maximum principle for superharmonic functions implies that $u-u_k \ge 0$ 
almost everywhere in $B$. Thus Fatou's Lemma gives
\[
 u(x) \ge \liminf_{k \rightarrow \infty} u_k(x) \ge \int_B G^x_B(y)f(y)~\!dy,
\]
for almost every $x \in B$, as claimed.

We will use (\ref{eq:sysIne}) to estimate the quantities
$$
{\cal U_R}=\int_{B_R}|x|^b|u|^{q-1}~\!dx~,\quad {\cal V_R}=\int_{B_R}|x|^a|v|^{p-1}~\!dx
$$
for any $R>0$ small enough (recall that ${\cal U_R}, {\cal V_R}$ are finite).
For $|x|<1/2$ and $y\in B$ we have the uniform lower bound
$$
G^x(y)\ge c_n\left[ \frac{1}{(|x|+|y|)^{n-2}} - M\right]~,\quad M=\max_{|x|\le \frac{1}{2}~,~ |y|<1} h^x(y).
$$
Therefore, if $R_0$ is small enough and $R\in(0,R_0)$, then $G^x_B(y) \ge cR^{2-n}$
for  $x, y\in B_R$.  Using \eqref{eq:sysIne}, we infer that
$$
u(x)\ge c R^{2-n}\int_{B_R}  | y |^av(y)^{p-1}~\!dy~,\quad v(x)\ge c R^{2-n}\int_{B_R}  | y |^bu(y)^{q-1}~\!dy
$$
for almost every $x\in B_R$, and hence
$$
{\cal U}_R\ge   cR^{(2-n)(q-1)+b+n}~\!{\cal V}_R^{q-1}~,\quad
\mathcal V_R\ge c R^{(2-n)(p-1) + a+n}~\! \mathcal U_R^{p-1}~\!.
$$
With simple computations we arrive at
\begin{gather}
\label{eq:UR}
{\cal U}_R^{(p-1)(q-1)-1}\le c~\!R^{-(q-1)p\left[\frac{a+n}{p}+\frac{b+n}{p(q-1)}-(n-2)]\right]}~\!,\\
\label{eq:VR}
{\cal V}_R^{(p-1)(q-1)-1}\le c~\!R^{-(p-1)q\left[\frac{a+n}{q(p-1)}+\frac{b+n}{q}-(n-2)]\right]}~\!.
\end{gather}
Now we distinguish three cases, depending whether ({C}), (H) or (AC) is satisfied.

\bigskip
\noindent
{\bf Coercive case ({C}).} We have that $\theta:=1-(p-1)(q-1)>0$, and thus (\ref{eq:VR})
gives
$$
R^{\frac{(p-1)q}{\theta}\left[\frac{a+n}{q(p-1)}+\frac{b+n}{q}-(n-2)]\right]}\le c{\cal V}_R.
$$
But clearly ${\cal V}_R\to 0$ as $R\to 0$. Thus
$$
\frac{a+n}{q(p-1)}+\frac{b+n}{q}-(n-2)>0.
$$
A similar argument and (\ref{eq:UR}) lead to the conclusion that $(a,b)\in E_{p,q}$.

\bigskip

\noindent{\bf Homogeneous case (H).} We have that $q=p'$, and therefore \eqref{eq:VR}  gives
$$
1\le c~\!R^{-p\left[\frac{a}{p}+\frac{b}{p'}+2\right]}.
$$ 
Hence $\frac{a}{p}+\frac{b}{p'}+2\ge 0$, that is, $(a,b)\in E_{p, p'}$.

\bigskip

\noindent
{\bf Anticoercive case (AC).} 
As $u, v$ are positive and superharmonic, they
are uniformly bounded from
below on any  ball $B_{R_0}\subset\Omega$ about the origin. Hence, for $R\in(0,R_0]$ we get
$$
{\cal V}_R=\int\limits_{B_R}|x|^{a}v^{p-1}~\!dx\ge c R^{n+a}~,\quad
{\cal U}_R=\int\limits_{B_R}|x|^{b}u^{q-1}~\!dx\ge c R^{n+b}~\!,
$$
that compared with (\ref{eq:UR}), (\ref{eq:VR}) give
$$
c\le ~\!R^{-(q-1)p\left[\frac{a}{p}+\frac{b}{p'}+2\right]}~,\qquad
c\le ~\!R^{-(p-1)q\left[\frac{a}{q'}+\frac{b}{q}+2\right]}~\!,
$$
as  $(p-1)(q-1)-1>0$. We immediately infer that
\begin{equation}
\label{eq:almost}
\frac{a}{q'}+\frac{b}{q}+2\ge 0~,\quad \frac{a}{p}+\frac{b}{p'}+2 \ge 0~\!.
\end{equation}

It remains to  prove that  strict inequalities hold in \eqref{eq:almost}.
We argue by contradiction. Assume for instance that
\begin{equation}
\label{eq:contra}
\frac{a}{q'}+\frac{b}{q}+2 = 0.
\end{equation}
Then the second inequality in \eqref{eq:almost} and $(p-1)(q-1)-1>0$
imply that $a\le -2\le b$. In addition, \eqref{eq:VR} becomes
\begin{equation}
\label{eq:VR2}
{\cal V}_R\le c~\!R^{n+a}~\!.
\end{equation}
We first consider the case
$$
-n<a<-2~,\qquad b=-2q-a(q-1)>-2.
$$
Since $v$ is bounded from below on a small ball $B_{2\sqrt{R}}$, then $-\Delta u\ge c|x|^a$ on $B_{2\sqrt{R}}$.
Thus, by the maximum principle, 
$$
u(x)\ge c\big(|x|^{a+2}-(2\sqrt{R})^{a+2}\big)\qquad\text{on ~$B_{2\sqrt{R}}$.}
$$
In particular, $u(x)\ge c~\!|x|^{a+2}$  on ~$B_{\sqrt{R}}$, so that
$$
-\Delta v\ge c~\!|x|^{-2}\qquad\text{on ~$B_{\sqrt{R}}$,}
$$
as $b+(a+2)(q-1)=-2$. Again by the maximum principle, $v(x)\ge c\log(\sqrt{R}/|x|)$ on $B_{\sqrt{R}}$,
and in particular
$$
v(x)\ge c\big|\log|x|\big|\qquad\text{on ~$B_{R}$.}
$$
We infer that
$$
\mathcal V_R\ge c\int_{B_R}|x|^a\big|\log|x|\big|^{p-1}~\!dx= O\Big(R^{n+a}|\log R|^{p-1}\Big),
$$
in contradiction with \eqref{eq:VR2}. To exclude the case $a=-2$, notice that in this case
$b=-2$ by \eqref{eq:contra},
hence $-\Delta v\ge c|x|^{-2}$ on $\Omega$. Conclude as before.
\QED

\subsection*{Acknowledgements}
The authors are pleased to thank Enzo Mitidieri and Lorenzo D'Ambrosio for useful remarks on how to 
improve the presentation of this paper.

\appendix

\section{\!\!\!\!\!\!ppendix: power-type solutions}
\label{S:elementary}

In \cite{BV}, Bidaut-Veron computed the positive ``power-type'' solutions $u, v$ to
\begin{equation}
\label{eq:power}
\begin{cases}
-\Delta u= \lambda_1|x|^av^{p-1}\\
-\Delta v= \lambda_2|x|^bu^{q-1}
\end{cases} 
\end{equation}
on $\R^n\setminus\{0\}$.  Actually we are just interested in finding
the set ${\mathcal Q_{p,q}}$
of parameters $a,b$,  for which (\ref{eq:HLE}) admits power-type solutions.

Let us start with a few remarks about the {\em Kelvin transform} 
defined in (\ref{eq:kelvin}). A simple computation shows that
$\mathcal K$ maps distributional
solutions $u,v\in L^1_\mathrm{loc}(\Omega\setminus\{0\})$ of \eqref{eq:HLE} 
into distributional solutions $\mathcal K u,\mathcal Kv\in L^1_\mathrm{loc}(\hat\Omega)$ to
$(\mathcal P_{\kappa(a,b)})$, where $\hat\Omega$ is the reflection of $\Omega$ with respect to the unitary sphere
and
$$
\kappa(a,b)=\left(-a-2n+p(n-2),-b-2n+q(n-2)\right)~,\qquad \kappa:\R^2\to \R^2~\!.
$$
Trivially, $\mathcal K$ maps power-type solutions into power-type solutions, that is, ${\mathcal Q_{p,q}}$
is invariant under the action of $\kappa$. Next we notice that $\kappa$ is a central inversion with respect to its fixed point $F$:
$$
\kappa(a,b)=2F-(a,b)~,\qquad F=\left(p\frac{n-2}{2}-n,q\frac{n-2}{2}-n\right)~\!.
$$

If the pair $u(x)=|x|^\alpha$, $v(x)=|x|^\beta$ solves \eqref{eq:power} with respect 
to some $\lambda_1, \lambda_2>0$, then clearly $\alpha, \beta$ have
to satisfy 
\begin{equation}
\label{eq:sys}
\begin{cases}
(q-1)\alpha-\beta=-b-2\\
-\alpha+(p-1)\beta=-a-2~\!.
\end{cases}
\end{equation}
In the non-homogeneous cases (AC) and ({C}), we have that the system \eqref{eq:sys} admits the unique solution
$$
\alpha=-\frac{\frac{a}{p}+\frac{b}{p'}+2}{(p-1)(q-1)-1}~\!p~,\quad 
\beta=-\frac{\frac{a}{q'}+\frac{b}{q}+2}{(p-1)(q-1)-1}~\!q~\!.
$$
The corresponding pair $u_\alpha, v_\beta$ solves \eqref{eq:HLE} with $\lambda_1,\lambda_2$
given, up to positive multipliers, by
\begin{eqnarray*}
\lambda_1&=&-\left(\frac{a}{p}+\frac{b}{p'}+2\right)\left(\frac{a+n}{q(p-1)}+\frac{b+n}{q}-(n-2)\right)\\
~\\
\lambda_2&=&-\left(\frac{a}{q'}+\frac{b}{q}+2\right)\left(\frac{a+n}{p}+\frac{b+n}{p(q-1)}-(n-2)\right).
\end{eqnarray*}
In conclusion, nontrivial and positive power-type solutions to \eqref{eq:power} exist 
if and only if the couple of exponents $(a, b)$ belongs to the open parallelogram ${\mathcal Q_{p,q}}$
whose vertices are
\begin{eqnarray*}
X=~(-n,q(n-2)-n),&&\quad X'=\kappa(X)=( p (n-2) - n, -n ),\\ V=(-2,-2),&&\quad V'=\kappa(V).
\end{eqnarray*}
More explicitly, if (AC) holds we have that 
$$
{\mathcal Q_{p,q}}=\left\{(a,b)\in \R^2~~:~~\begin{array}{l}
\min\left\{\displaystyle{\frac{a}{p}+\frac{b}{p'}+2},\displaystyle{\frac{a}{q'}+\frac{b}{q}+2}\right\} >  0\\
~\\
\max\left\{\displaystyle{\frac{a+n}{p}+\frac{b+n}{p(q-1)}}, \displaystyle{\frac{a+n}{q(p-1)}+\frac{b+n}{q}}\right\} < n-2
\end{array}\right\},
$$
while in the coercive case ({C}) we find
$$
{\mathcal Q_{p,q}}=\left\{(a,b)\in \R^2~~:~~\begin{array}{l}
\max\left\{\displaystyle{\frac{a}{p}+\frac{b}{p'}+2},\displaystyle{\frac{a}{q'}+\frac{b}{q}+2}\right\} <  0\\
~\\
\min\left\{\displaystyle{\frac{a+n}{p}+\frac{b+n}{p(q-1)}}, \displaystyle{\frac{a+n}{q(p-1)}+\frac{b+n}{q}}\right\} > n-2
\end{array}\right\}.
$$
Points in the boundary of $\mathcal Q_{p,q}$ correspond to trivial solutions to \eqref{eq:power} in the sense 
of Bidaut-Veron~\cite{BV}, that is, at least one of the components is harmonic on $\R^n\setminus\{0\}$. 

The coordinates of the vertices $X, X'$ satisfy
\begin{equation*}\tag{\em CL}
\label{eq:CL}
\frac{a+n}{p}+\frac{b+n}{q}=n-2~\!.
\end{equation*}
The remaining vertices $V$ and $V'$ lie on opposite sides of line in the $a,b$ plane given by (\ref{eq:CL}).
More precisely, $V$ is below (\ref{eq:CL}) in the anticoercive case (AC), while $V$ is above (\ref{eq:CL}) if ({C}) holds.

\bigskip
In the homogenous case (H), the line \eqref{eq:CL} becomes
\begin{equation*}\tag{$\textit{CL}_{\text H}$}
\label{eq:NCL}
\frac{a}{p}+\frac{b}{p'}+2=0~\!
\end{equation*}
and with simple calculations we find that
${\mathcal Q_{p,q}}$ collapses into
$$
{\mathcal Q_{p,p'}}=\left\{(a,b)\in \R^2~:~a,b>-n~,~~\displaystyle{\frac{a}{p}+\frac{b}{p'}+2}=0\right\}~\!,
$$
that is the open segment of endpoints $X, X'$. 

In the next pictures we represent the  
set $\mathcal Q_{p,q}$ in the $(a,b)$ plane.
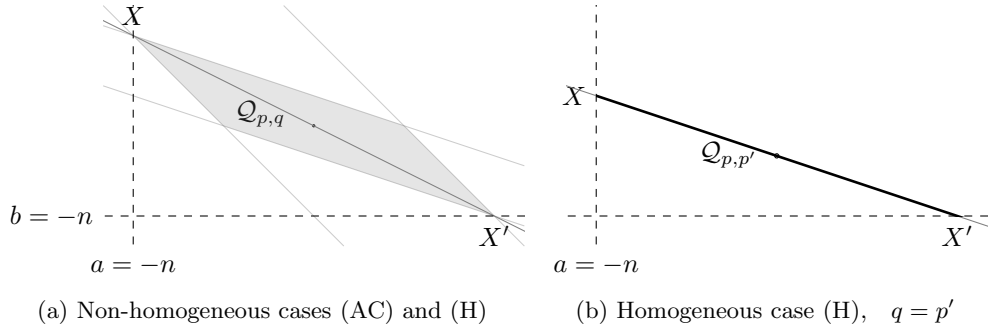
\begin{figure}[H]
\label{F:Q2}
\centering
\subfloat[Non-homogeneous cases (AC) and (H)]{
\begin{tikzpicture}[scale=.4]
\path (-6, -5) node[left]{\small $b=-n$};
\path (-5, -6) node[below]{\small $a=-n$};
\draw (1, -2) circle(1pt);
\path (0.4, -1.6) node[left]{\small ${\mathcal Q_{p,q}}$};
\path (-5, 1) node[above]{\small $X$};
\path (7, -5) node[below]{\small $X'$};
\clip (-6, -6) rectangle (8, 2);
\draw[color=gray] (-7, 2)--(11, -7);
\draw[color=gray, opacity=.4] (-6, -0.67)--(8, -5.33);
\draw[color=gray, opacity=.4] (-6, 1.33)--(8, -3.33);
\draw[color=gray, opacity=.4] (2, -6)--(-6, 2);
\draw[color=gray, opacity=.4] (8, -6)--(0, 2);
\draw[dashed] (-7, -5)--(9, -5);
\draw[dashed] (-5, -7)--(-5, 3);
\clip (-6, -0.67)--(8, -5.33)--(8, 2)--(-6, 2)--cycle;
\clip (-6, 1.33)--(8, -3.33)--(8, -6)--(-6, -6)--cycle;
\clip (2, -6)--(-6, 2)--(8, 2)--(8, -6)--cycle;
\clip (8, -6)--(0, 2)--(-6, 2)--(-6, -6)--cycle;
\fill[color=gray, opacity=.2] (-6, -6) rectangle (8, 2);
\end{tikzpicture}
}
\subfloat[Homogeneous case (H),~~ $q=p'$]{
 \begin{tikzpicture}[scale=.4]
\path (-5, -6) node[below]{\small $a=-n$};
\draw (1, -3) circle(2pt);
\path (-0.6, -2.2) node[below]{\small ${\mathcal Q_{p,p'}}$};
\path (-5, -1) node[left]{\small $X$};
\path (7, -5) node[below]{\small $X'$};
\clip (-6, -6) rectangle (8, 2);
\draw[dashed] (-7, -5)--(9, -5);
\draw[dashed] (-5, -7)--(-5, 3);
\draw[color=gray] (22, -10)--(-38, 10);
\clip (-5, -5) rectangle (8, 2);
\draw[line width=1pt] (22, -10)--(-38, 10);
\end{tikzpicture}
}
\end{figure}

\noindent Notice that, in any case, $\mathcal Q_{p,q}= E_{p,q}\cap \kappa(E_{p,q})$. 

In the next pictures we summarize our existence/nonexistence results.
We have existence of weak solutions on bounded neighborhoods of the origin
in the light gray zone. Power-type solutions correspond to the darker area. 
The Brezis-Cabr\'e
nonexistence result for the inequality (\ref{eq:BC}) is related to the vertex $V=(-2,-2)$ in Figure $({c})$
(with $p=q>2$).

\begin{figure}[H]
\label{F:E}
\centering
\subfloat[Anticoercive case~~ $\frac{1}{p}+\frac{1}{q}< 1$]{
\begin{tikzpicture}[scale=.4]
\path (-5, 1) node[left]{\small $X$};
\path (7, -5) node[below]{\small $X'$};
\clip (-6, -6) rectangle (8, 2);
 \draw[color=gray] (-7, 2)--(11, -7);
\draw[color=gray] (-2, -2)--(8, -5.33);
\draw[color=gray] (-2, -2)--(-6, 2);
\path (-2, -2) node[below]{\footnotesize $V$};
\draw (-2, -2) circle(1pt);
\draw[dashed] (-7, -5)--(9, -5);
\draw[dashed] (-5, -7)--(-5, 3);
\clip (-6, -0.67)--(8, -5.33)--(8, 2)--(-6, 2)--cycle;
\clip (2, -6)--(-6, 2)--(8, 2)--(8, -6)--cycle;
\clip (-5, -5) rectangle (8, 2);
\fill[color=gray, opacity=.2] (-6, -6) rectangle (8, 2);
\path (4, 0) node{\small $E_{p, q}$};
\clip (8, -6)--(0, 2)--(-6, 2)--(-6, -6)--cycle;
\clip (-6, 1.33)--(8, -3.33)--(8, -6)--(-6, -6)--cycle;
\fill[color=gray, opacity=.4] (-6, -6) rectangle (8, 2);
\end{tikzpicture}
}\qquad
\subfloat[Coercive case~~ $\frac{1}{p}+\frac{1}{q}> 1$]{
\begin{tikzpicture}[scale=.4]
\path (-5, 1) node[left]{\small $X$};
\path (7, -5) node[below]{\small $X'$};
\clip (-6, -6) rectangle (8, 2);
 \draw[color=gray] (-7, 2)--(11, -7);
\draw[color=gray] (-2, -2)--(8, -5.33);
\draw[color=gray] (-2, -2)--(-6, 2);
\path (4, -2) node[above]{\footnotesize $V$};
\draw[dashed] (-7, -5)--(9, -5);
\draw[dashed] (-5, -7)--(-5, 3);
\clip (-6, -0.67)--(8, -5.33)--(8, 2)--(-6, 2)--cycle;
\clip (2, -6)--(-6, 2)--(8, 2)--(8, -6)--cycle;
\clip (-5, -5) rectangle (8, 2);
\fill[color=gray, opacity=.2] (-6, -6) rectangle (8, 2);
\path (4, 0) node{\small $E_{p, q}$};
\clip (8, -6)--(0, 2)--(-6, 2)--(-6, -6)--cycle;
\clip (-6, 1.33)--(8, -3.33)--(8, -6)--(-6, -6)--cycle;
\fill[color=gray, opacity=.4] (-6, -6) rectangle (8, 2);
\end{tikzpicture}
}
\end{figure}

\end{document}